\documentclass[conference]{IEEEtran}
\IEEEoverridecommandlockouts
\usepackage{cite}
\usepackage{amsmath,amssymb,amsfonts}
\usepackage{algorithmic}
\usepackage{graphicx}
\usepackage{textcomp}
\usepackage{xcolor}
\def\BibTeX{{\rm B\kern-.05em{\sc i\kern-.025em b}\kern-.08em
    T\kern-.1667em\lower.7ex\hbox{E}\kern-.125emX}}

\usepackage{amsmath,amsfonts,amssymb,amsthm}
\theoremstyle{definition}

\theoremstyle{plain}

\newtheorem{lemma}{Lemma}
\newtheorem{proposition}{Proposition}

\usepackage{mathrsfs}
\usepackage{booktabs}
\usepackage{multirow}
\usepackage{array}
\usepackage[ruled]{algorithm}
\usepackage{hyperref}
\hypersetup{colorlinks,%
	citecolor={red}, 
	urlcolor={blue},
	linkcolor={blue},
	breaklinks={true},
	pagebackref={true},
	hyperindex={true},
}
\usepackage{url}
\usepackage{pgf,tikz}
\usetikzlibrary{arrows,shapes,positioning}
\usepackage{circuitikz}
\usepackage{subfig}
\usepackage{setspace}
\usepackage{url}
\floatname{algorithm}{Model}

\usepackage{enumitem}
\newlist{abbrv}{itemize}{1}
\setlist[abbrv,1]{label=,labelwidth=0.9in,align=parleft,noitemsep,leftmargin=!}

\newcommand{\R}{\mathbb{R}}
\newcommand{\C}{\mathbb{C}}

\newcommand{\Herm}{\mathbb{H}}

\newcommand{\cplx}[1]{\mathrm{#1}}
\newcommand{\rv}[1]{\boldsymbol{#1}}
\newcommand{\cv}[1]{\boldsymbol{\mathrm{#1}}}
\newcommand{\ub}[1]{\overline{#1}}
\newcommand{\lb}[1]{\underline{#1}}
\newcommand{\Node}{\mathcal{N}}
\newcommand{\Gen}{\mathcal{G}}
\newcommand{\Branch}{\mathcal{L}}
\newcommand{\Reseau}{\mathscr{P}}

\newcommand{\jc}{\cplx{j}}
\newcommand{\pl}{\lb{p}}
\newcommand{\ql}{\lb{q}}
\newcommand{\pu}{\ub{p}}
\newcommand{\qu}{\ub{q}}
\newcommand{\vl}{\lb{v}}
\newcommand{\vu}{\ub{v}}
\DeclareMathOperator{\subj}{s.t.}
\DeclareMathOperator{\rank}{rank}

\DeclareMathOperator{\re}{Re}
\DeclareMathOperator{\im}{Im}

\begin{document}

\title{CONICOPF: Conic Relaxations for AC Optimal Power Flow Computations
\thanks{This research was supported by the NSERC-Hydro-Qu\'ebec-Schneider Electric Industrial Research Chair.}
}

\author{\IEEEauthorblockN{Christian Bingane}
\IEEEauthorblockA{
\textit{Polytechnique Montr\'eal}\\
Montr\'eal, Qu\'ebec, Canada H3C~3A7 \\
\url{christian.bingane@polymtl.ca}}
\and
\IEEEauthorblockN{Miguel~F. Anjos}
\IEEEauthorblockA{
\textit{University of Edinburgh}\\
Edinburgh, Scotland, UK EH9~3FD \\
\url{anjos@stanfordalumni.org}}
\and
\IEEEauthorblockN{S\'ebastien {Le~Digabel}}
\IEEEauthorblockA{
\textit{Polytechnique Montr\'eal}\\
Montr\'eal, Qu\'ebec, Canada H3C~3A7 \\
\url{sebastien.le.digabel@gerad.ca}}
}
\maketitle

\begin{abstract}
Computational speed and global optimality are key needs for practical algorithms for the optimal power flow problem. 
Two convex relaxations offer a favorable trade-off between the standard second-order cone and the standard semidefinite relaxations for large-scale meshed networks in terms of optimality gap and computation time: the tight-and-cheap relaxation (TCR) and the quadratic convex relaxation (QCR). 
We compare these relaxations on 60 PGLib-OPF test cases with up to 1,354 buses  under three operating conditions 
and show that TCR dominates QCR on all 20 typical (TYP) test cases, on 18 out of 20 active power increase (API), and 12 out of 20 small angle difference (SAD). Selected state-of-the-art conic relaxations are implemented in the new MATLAB-based package CONICOPF available on GitHub.
\end{abstract}

\begin{IEEEkeywords}
Power systems, optimal power flow, convex optimization, conic relaxations.
\end{IEEEkeywords}

\section*{Nomenclature}
\subsection{Sets}
\begin{abbrv}
	\item[$\R$/$\C$] Set of real/complex numbers,
	\item[$\Herm^n$] Set of $n\times n$ Hermitian matrices,
	\item[$\jc$] Imaginary unit,
	\item[$a$/$\cplx{a}$] Real/complex number,
	\item[$\rv{a}$/$\cv{a}$] Real/complex vector,
	\item[$A$/$\cplx{A}$]  Real/complex matrix.
\end{abbrv}
\subsection{Operators}
\begin{abbrv}
	\item[$\re (\cdot)$/$\im (\cdot)$]  Real/imaginary part operator,
	\item[$(\cdot)^*$] Conjugate operator,
	\item[$\left|\cdot\right|$]  Magnitude or cardinality set operator,
	\item[$\angle (\cdot)$]  Phase operator,
	\item[$(\cdot)^H$]  Conjugate transpose operator,
	\item[$\rank (\cdot)$]  Rank operator.
\end{abbrv}
\subsection{Input data}
\begin{abbrv}
	\item[$\Reseau = (\Node, \Branch)$] Power network,
	\item[$\Node$] Set of buses,
	\item[$\Gen = \bigcup_{k\in\Node} \Gen_k$] Set of generators,
	\item[$\Gen_k$] Set of generators connected to bus~$k$,
	\item[$\Branch$] Set of branches,
	\item[$p_{Dk}$/$q_{Dk}$] Active/reactive power demand at bus~$k$,
	\item[$g_k'$/$b_k'$] Conductance/susceptance of shunt at bus~$k$, 
	\item[$c_{g2}, c_{g1}, c_{g0}$] Generation cost coefficients of generator~$g$,
	\item[$\cplx{y}_\ell^{-1} = r_{\ell} + \jc x_{\ell}$] Series impedance of branch~$\ell$,
	\item[$b'_\ell$] Total shunt susceptance of branch~$\ell$,
	\item[$\cplx{t}_{\ell}$] Turns ratio of branch~$\ell$.
\end{abbrv}
\subsection{Variables}
\begin{abbrv}
	\item[$p_{Gg}$/$q_{Gg}$] Active/reactive power output by generator~$g$,
	\item[$\cplx{v}_{k}$] Complex (phasor) voltage at bus~$k$,
	\item[$p_{f\ell}$/$q_{f\ell}$] Active/reactive power flow injected along branch~$\ell$ by its \emph{from} end,
	\item[$p_{t\ell}$/$q_{t\ell}$] Active/reactive power flow injected along branch~$\ell$ by its \emph{to} end.
\end{abbrv}

\section{Introduction}
\IEEEPARstart{T}{he} optimal power flow (OPF) problem seeks to find an optimal network operating point subject to power flow equations and other operational constraints~\cite{ref2}. The continuous classical version with AC power flow equations, which is nonconvex and NP-hard~\cite{lehmann2016}, is generally called AC optimal power flow (ACOPF) problem.

In recent years, convex relaxations of the ACOPF problem, such as the second-order cone relaxation (SOCR)~\cite{jabr1} and the semidefinite relaxation (SDR)~\cite{bai}, have attracted a significant interest for several reasons. First, they provide a bound on the global optimal value of ACOPF and may lead to global optimality. Second, if one of them is infeasible, then the ACOPF problem is infeasible. Convex relaxations do not replace nonlinear (local) solvers but rather complement them because they measure the quality of the solution obtained~\cite{capitanescu2016}.

For general meshed networks, SDR is stronger than SOCR but requires heavier computation. Therefore, the chordal relaxation (CHR) was proposed in~\cite{jabr4} in order to exploit the fact that power networks are not densely connected, thus reducing data storage and increasing computation speed. However, even CHR remains expensive to solve compared to SOCR for large-scale power systems. On the other hand, for radial networks, SOCR is tantamount to SDR. In this case, one would normally solve the first one rather than the second one due to the difference in computation time.
Recent surveys on these relaxations can be found in~\cite{molzahn2019,zohrizadeh2020}.

According to~\cite{frank1}, high computational speed is a key need for practical OPF algorithms, especially in real-time applications and when dealing with large-scale power systems. In fact, in real-time applications, an OPF problem is run every few minutes to update device and resource settings in response to the constantly changing conditions of power systems~\cite{taylor2015}. This need motivated the tight-and-cheap relaxation (TCR)~\cite{bingane2018,bingane2018tight,bingane2019phd} that offers a favorable trade-off between SOCR and SDR for large-scale meshed instances of ACOPF in terms of optimality gap and computation time. Indeed, TCR was proven to be stronger than SOCR and nearly as tight as SDR. Moreover, computational experiments on MATPOWER test cases with up to 6,515 buses showed that solving TCR for large-scale instances is much less expensive than solving CHR.

In the literature, there is another convex relaxation that offers the same advantages as TCR: the quadratic convex relaxation (QCR)~\cite{hijazi2014}. A theoretical and computational study about this relaxation can be found in~\cite{coffrin2016a}. It was shown that QCR is neither weaker nor stronger than SDR, but it is computationally cheaper to solve. 

In this paper, we compare TCR and QCR on 60 PGLib-OPF test cases, ranging from 3 to 1,354 buses and consisting of 20 different networks under three operating conditions: typical (TYP), active power increase (API), and small angle difference (SAD). The optimality gaps show that TCR on average dominates QCR. Moreover, under TYP operating conditions, TCR is stronger than QCR. To the best of our knowledge, it is the first time that such a comparison is made between both relaxations on PGLib-OPF instances, which are more challenging than MATPOWER's. 
We also introduce CONICOPF, a MATLAB-based package available on GitHub~\cite{conicopf} with selected state-of-the-art conic relaxations.

This paper is organized as follows. In Section~\ref{sec3:tcr}, we recall TCR and QCR. We report in Section~\ref{sec3:results} computational results on the comparison between TCR and QCR. Section~\ref{sec3:conclusion} concludes the paper.

\section{ACOPF: Convexification}\label{sec3:tcr}
Consider a typical power network $\Reseau = (\Node, \Branch)$, where $\Node = \{1,2,\ldots,n\}$ and $\Branch \subseteq \Node \times \Node$ denote respectively the set of buses and the set of branches (transmission lines, transformers and phase shifters). Each branch~$\ell \in \Branch$ has a \emph{from} end  $k$ (on the \emph{tap side}) and a \emph{to} end $m$  as modeled in~\cite{matpower}. We note $\ell = (k,m)$.  The ACOPF problem is given as:
\begin{subequations}\label{eq3:acopf}
	\begin{equation}\label{eq3:objfun}
		\min \sum_{g\in\Gen} c_{g2} p_{Gg}^{2} + c_{g1} p_{Gg} + c_{g0}
	\end{equation}
	over variables $\rv{p}_{G}, \rv{q}_{G} \in \R^{|\Gen|}$, $\rv{p}_{f},\rv{q}_{f},\rv{p}_{t},\rv{q}_{t} \in \R^{|\Branch|}$, and $\cv{v} \in~\C^{|\Node|}$, subject to
	\begin{itemize}
		\item Power balance equations:
		\begin{align}
			&\sum_{g\in \Gen_k} p_{Gg} - p_{Dk} - g_k'\left|\cplx{v}_k\right|^2 = \nonumber\\
			&\sum_{\ell=(k,m)\in\Branch} p_{f\ell} + \sum_{\ell=(m,k)\in\Branch} p_{t\ell}\;\forall k\in\Node,\label{eq3:kclp}\\
			&\sum_{g\in \Gen_k} q_{Gg} - q_{Dk} + b_k'\left|\cplx{v}_k\right|^2 = \nonumber\\ &\sum_{\ell=(k,m)\in\Branch} q_{f\ell} + \sum_{\ell=(m,k)\in\Branch} q_{t\ell}\;\forall k\in\Node,\label{eq3:kclq}
		\end{align}
		\item Line flow equations:
		\begin{align}
			\frac{\cplx{v}_k}{\cplx{t}_{\ell}} &\left[\left(\jc \frac{b'_{\ell}}{2} + \cplx{y}_\ell \right) \frac{\cplx{v}_k}{\cplx{t}_{\ell}} - \cplx{y}_\ell \cplx{v}_m\right]^*\nonumber\\
			&= p_{f\ell} + \jc q_{f\ell}\;\forall \ell = (k,m)\in\Branch,\label{eq3:sf}\\
			\cplx{v}_m &\left[- \cplx{y}_\ell \frac{\cplx{v}_k}{\cplx{t}_{\ell}} + \left(\jc \frac{b'_{\ell}}{2} + \cplx{y}_\ell \right) \cplx{v}_m\right]^*\nonumber\\
			&= p_{t\ell} + \jc q_{t\ell}\;\forall \ell = (k,m)\in\Branch,\label{eq3:st}
		\end{align}
		\item Generator power capacities:
		\begin{equation}\label{eq3:genlim}
			\pl_{Gg}\le p_{Gg} \le \pu_{Gg},\, \ql_{Gg}\le q_{Gg} \le \qu_{Gg}\;\forall g\in\Gen,
		\end{equation}
		\item Line thermal limits:
		\begin{equation}\label{eq3:linelim}
			|p_{f\ell} + \jc q_{f\ell}| \le \overline{s}_{\ell},\, |p_{t\ell} + \jc q_{t\ell}| \le \ub{s}_{\ell}\;\forall \ell \in\Branch,
		\end{equation}
		\item Voltage magnitude limits:
		\begin{equation}\label{eq3:magbus}
			\vl_k\le \left|\cplx{v}_k\right| \le \vu_k\;\forall k\in\Node,
		\end{equation}
		\item Phase angle difference limits:
		\begin{equation}\label{eq3:phaselim}
			|\angle \cplx{v}_k - \angle \cplx{v}_m| \le \ub{\delta}_{\ell}\;\forall \ell = (k,m) \in\Branch,
		\end{equation}
		\item Reference bus constraint:
		\begin{equation}\label{eq3:slack}
			\angle \cplx{v}_1 = 0.
		\end{equation}
	\end{itemize}
\end{subequations}
The objective function~\eqref{eq3:objfun} is the cost of conventional generation commonly used in the literature. Constraints~\eqref{eq3:kclp}--\eqref{eq3:st} are derived from Kirchhoff's laws and represent power flows in the network. Constraint~\eqref{eq3:slack} specifies bus $k=1$ as the reference bus. We assume that $\vl_k > 0$ for all $k\in \Node$ in~\eqref{eq3:magbus}, that $0 < \ub{\delta}_\ell < \pi/2$ for all $\ell \in \Branch$ in~\eqref{eq3:phaselim}, and that the generation cost $c_{g2} p_{Gg}^2 + c_{g1} p_{Gg} + c_{g0}$ is a convex function for all $g\in\Gen$.

Problem~\eqref{eq3:acopf} is highly nonconvex and NP-hard~\cite{lehmann2016} due to the nonconvex constraints~\eqref{eq3:sf}--\eqref{eq3:st}.

\subsection{Tight-and-cheap relaxation}
With $\cplx{V} := \cv{v} \cv{v}^{H}$, ACOPF~\eqref{eq3:acopf} can be reformulated as:
\begin{subequations}\label{eq3:acopfV}
	\begin{align}
		\min\quad &\eqref{eq3:objfun} \nonumber\\
		\subj\quad &\eqref{eq3:genlim},~\eqref{eq3:linelim},~\eqref{eq3:slack}, \nonumber\\
		&\sum_{g\in \Gen_k} p_{Gg} - p_{Dk} - g_k'\cplx{V}_{kk} = \nonumber\\
		&\sum_{\ell=(k,m)\in\Branch} p_{f\ell} + \sum_{\ell=(m,k)\in\Branch} p_{t\ell}\;\forall k\in\Node,\label{eq3:kclpV}\\
		&\sum_{g\in \Gen_k} q_{Gg} - q_{Dk} + b_k'\cplx{V}_{kk} = \nonumber\\
		&\sum_{\ell=(k,m)\in\Branch} q_{f\ell} + \sum_{\ell=(m,k)\in\Branch} q_{t\ell}\;\forall k\in\Node, \label{eq3:kclqV}\\
		&\frac{1}{|\cplx{t}_{\ell}|^2} \left(-\jc \frac{b'_{\ell}}{2} + \cplx{y}_\ell^* \right) \cplx{V}_{kk} - \frac{\cplx{y}_\ell^*}{\cplx{t}_{\ell}} \cplx{V}_{km} \nonumber\\
		&= p_{f\ell} + \jc q_{f\ell}\;\forall \ell = (k,m)\in\Branch,\label{eq3:sfV}\\
		&-\frac{\cplx{y}_\ell^*}{\cplx{t}_{\ell}^*} \cplx{V}_{km}^* + \left(-\jc \frac{b'_{\ell}}{2} + \cplx{y}_\ell^* \right) \cplx{V}_{mm} \nonumber\\
		&= p_{t\ell} + \jc q_{t\ell}\;\forall \ell = (k,m)\in\Branch, \label{eq3:stV}\\
		&\vl_k^2\le \cplx{V}_{kk} \le \vu_k^2\;\forall k\in\Node, \label{eq3:Vkk}\\
		&|\im (\cplx{V}_{km})| \le \re (\cplx{V}_{km})\tan \ub{\delta}_\ell\; \forall \ell = (k,m)\in\Branch, \label{eq3:Vkm}\\
		&\cplx{V} = \cv{v} \cv{v}^{H}.\label{eq3:Vv}
	\end{align}
\end{subequations}
The nonconvexity of~\eqref{eq3:acopfV} is captured by~\eqref{eq3:Vv}. 
The standard semidefinite relaxation (SDR)~\cite{bai} is obtained using the equivalence $\cplx{V} = \cv{v} \cv{v}^{H}$ iff $\cplx{V} \succeq 0$ and $\rank (\cplx{V}) = 1$, and dropping the rank constraint. 
If we relax $\cplx{V} \succeq 0$ in SDR by $|\Branch|$ constraints of the form
\begin{equation}\label{eq3:socr}
	\cplx{V}_{\{k,m\}} :=
	\begin{bmatrix}
		\cplx{V}_{kk} & \cplx{V}_{km}\\
		\cplx{V}_{km}^* & \cplx{V}_{mm}
	\end{bmatrix}
	\succeq 0\; \forall (k,m) \in \Branch,\\
\end{equation}
we obtain the standard second-order cone relaxation (SOCR)~\cite{jabr1}, which is equivalent to SDR for radial networks.

SDR is very expensive to solve for large-scale instances, while SOCR is weaker than SDR for meshed networks. 
SDR can be solved more efficiently by exploiting the sparsity of a power network to replace the semidefinite constraint~$\cplx{V} \succeq 0$ by smaller 
semidefinite constraints based on a chordal extension of the network, thus obtaining the chordal relaxation (CHR) equivalent to SDR for meshed networks~\cite{jabr4}. In spite of the significant computational speed-up for large-scale instances obtained with CHR, this remains expensive compared to SOCR. This motivated the \emph{tight-and-cheap relaxation} (TCR)~\cite{bingane2018} obtained by combining semidefinite optimization with the reformulation-linearization technique (RLT).

The nonconvex constraint~\eqref{eq3:Vv} can be rewritten as
\begin{subequations}\label{eq3:Vvrectangular}
\begin{align}
	\cplx{V}_{kk} &= |\cplx{v}_{k}|^2 &\forall k \in \Node,\\
	\cplx{V}_{km} &= \cplx{v}_{k}\cplx{v}_{m}^{*} &\forall (k,m)\in\Branch.
\end{align}
\end{subequations}
To construct TCR in Model~\ref{model3:tcr}, we replace~\eqref{eq3:Vvrectangular} by the semidefinite constraint
\begin{subequations}\label{eq3:tcr}
	\begin{equation}
		\begin{bmatrix}
			1 & \cplx{v}_{k}^{*} & \cplx{v}_{m}^{*}\\
			\cplx{v}_{k} &\cplx{V}_{kk} & \cplx{V}_{km}\\
			\cplx{v}_{m} & \cplx{V}_{km}^{*} & \cplx{V}_{mm}
		\end{bmatrix} \succeq 0\; \forall \ell = (k,m)\in\Branch,\label{eq3:Vvtcr}
	\end{equation}
	and we add the following RLT constraints
	\begin{align}\label{eq3:rlt}
		\re (\cplx{v}_1) \ge \frac{\cplx{V}_{11} + \vl_1 \vu_1}{\vl_1 + \vu_1}, \quad
		\im (\cplx{v}_1) = 0,
	\end{align}
	corresponding to the reference bus~$k=1$.
\end{subequations}

\begin{algorithm}
	\caption{Tight-and-cheap relaxation (TCR)}
	\label{model3:tcr}
	\begin{algorithmic}
		\STATE Variables: $\rv{p}_G, \rv{q}_G \in \R^{|\Gen|}$, $\rv{p}_f, \rv{q}_f, \rv{p}_t, \rv{q}_t \in\R^{|\Branch|}$, $\cv{v}  \in\C^{|\Node|}$, and $\cplx{V}\in\Herm^{|\Node|}$.
		\STATE Minimize:~\eqref{eq3:objfun}
		\STATE Subject to:~\eqref{eq3:genlim},~\eqref{eq3:linelim},~\eqref{eq3:kclpV}--\eqref{eq3:Vkm},~\eqref{eq3:tcr}.
	\end{algorithmic}
\end{algorithm}

\subsection{Quadratic convex relaxation}
There is another convex relaxation offering the same advantages as TCR: the quadratic convex relaxation (QCR). This relaxation also provides an interesting alternative to SDR in that it is computationally cheaper and can be stronger than SDR when phase angle difference bounds in~\eqref{eq3:phaselim} are tight. While TCR preserves implicitly the links~\eqref{eq3:Vv} between $\cv{v}$ and~$\cplx{V}$ in rectangular form, QCR was introduced in~\cite{hijazi2014} to preserve these links in polar form.

With $v_k := |\cplx{v}_k|$ and $\theta_k := \angle \cplx{v}_k$ for all $k \in \Node$, the polar representation of the nonconvex constraint~\eqref{eq3:Vv} is given by
\begin{subequations}\label{eq3:Vvpolar}
	\begin{align}
		\cplx{V}_{kk} &= v_k^2 &\forall k \in \Node,\\
		\re (\cplx{V}_{km}) &= v_k v_m \cos (\theta_k - \theta_m) &\forall \ell = (k,m)\in \Branch,\\
		\im (\cplx{V}_{km}) &= v_k v_m \sin (\theta_k - \theta_m) &\forall \ell = (k,m)\in \Branch.
	\end{align}
\end{subequations}
Based on the following lemmas, QCR then formulates quadratic convex constraints to relax~\eqref{eq3:Vvpolar}.
\begin{lemma}
	Let $x$ be a real variable such that $\lb{x} \leq x \leq \ub{x}$,  where $\lb{x}< \ub{x}$. If $y = x^2$, then $x^2 \le y \le (\lb{x} + \ub{x})x - \lb{x}\ub{x}$.
\end{lemma}
\begin{lemma}
	For two real variables $x$, $y$ such that $\lb{x} \leq x \leq \ub{x}$, $\lb{y} \leq y \leq \ub{y}$,  where $\lb{x}< \ub{x}$, $\lb{y}< \ub{y}$, if $z = xy$ then $\max\{x\lb{y}+\lb{x}y-\lb{x}\lb{y}, x\ub{y}+\ub{x}y-\ub{x}\ub{y}\} \le z \le \min\{x\lb{y}+\ub{x}y-\ub{x}\lb{y}, x\ub{y}+\lb{x}y-\lb{x}\ub{y}\}$.
\end{lemma}
\begin{lemma}
	Let $t \in \R$ such that $|t| \le \ub{t}$, where $0 < \ub{t} < \pi/2$. If $x = \cos t$ and $y = \sin t$, then $\cos \ub{t} \le x \le 1 - (1-\cos \ub{t}) \frac{t^2}{\ub{t}^2}$ and $\left|y - t\cos \frac{\ub{t}}{2}\right| \le \sin \frac{\ub{t}}{2} - \frac{\ub{t}}{2} \cos \frac{\ub{t}}{2}$.
\end{lemma}

Now, let $w_\ell := v_k v_m$, $\tilde{c}_\ell := \cos (\theta_k - \theta_m)$ and $\tilde{s}_\ell := \sin (\theta_k - \theta_m)$ for all $\ell = (k,m) \in \Branch$. Assuming that $\lb{v}_k \lb{v}_m \le w_\ell \le \ub{v}_k \ub{v}_m$, $\cos \ub{\delta}_\ell \le \tilde{c}_\ell \le 1$, and $|\tilde{s}_\ell| \le \sin \ub{\delta}_\ell$ for all $\ell = (k,m) \in \Branch$, we construct QCR in Model~\ref{model3:qcr}. For a nonconvex constraint of the form $y = f(\rv{x})$, we note $y = \langle f(\rv{x})\rangle$ its convex relaxation.

\begin{algorithm}
	\caption{Quadratic convex relaxation (QCR)}
	\label{model3:qcr}
	\begin{algorithmic}
		\STATE Variables: $\rv{p}_G, \rv{q}_G \in \R^{|\Gen|}$, $
			\rv{p}_f, \rv{q}_f, \rv{p}_t, \rv{q}_t, \rv{w}, \tilde{\rv{c}}, \tilde{\rv{s}} \in\R^{|\Branch|}$, $	\rv{v}, \rv{\theta}  \in\R^{|\Node|}$, and $\cplx{V}  \in\Herm^{|\Node|}$.
		\STATE Minimize:~\eqref{eq3:objfun}
		\STATE Subject to:~\eqref{eq3:genlim},~\eqref{eq3:linelim},~\eqref{eq3:kclpV}--\eqref{eq3:stV},~\eqref{eq3:Vkm},~\eqref{eq3:socr}, and
		\begin{subequations}
			\begin{align}
				&\cplx{V}_{kk} = \langle v_k^2 \rangle &\forall k \in \Node,\\
				&\re (\cplx{V}_{km}) = \langle w_\ell \tilde{c}_\ell \rangle &\forall \ell = (k,m)\in \Branch,\\
				&\im (\cplx{V}_{km}) = \langle w_\ell \tilde{s}_\ell \rangle &\forall \ell = (k,m)\in \Branch,\\
				&w_\ell = \langle v_kv_m \rangle &\forall \ell = (k,m)\in \Branch,\\
				&\tilde{c}_\ell = \langle \cos (\theta_k - \theta_m) \rangle &\forall \ell = (k,m)\in \Branch,\\
				&\tilde{s}_\ell = \langle \sin (\theta_k - \theta_m) \rangle &\forall \ell = (k,m)\in \Branch,\\
				&\theta_1 = 0.
			\end{align}
		\end{subequations}
	\end{algorithmic}
\end{algorithm}

The computational results in~Table~\ref{table3:cost} show that 
\begin{proposition}
	TCR is neither weaker nor stronger than QCR.
\end{proposition}

\section{Computational results}\label{sec3:results}
\subsection{CONICOPF}
TCR, QCR, and others relaxations defined in~\cite{bingane2018,bingane2018tight,bingane2019phd} were implemented as a freely available MATLAB package: CONICOPF~\cite{conicopf}. 
CONICOPF provides M-files for solving relaxations of the classical ACOPF problem, the optimal reactive power dispatch (ORPD) problem, and the multi-period ACOPF problem. 
CONICOPF requires that MATPOWER~\cite{matpower} and CVX~\cite{cvx2} be installed and the input test case be in MATPOWER format.

\subsection{TCR vs QCR}
In this section, we compare TCR and QCR.
We tested Model~\ref{model3:tcr} and Model~\ref{model3:qcr} on test cases available from PGLib-OPF~19.05~\cite{pglibopf}, ranging from 3 to 1,354 buses and consisting of 20 different networks under typical (TYP), active power increase (API), and small angle difference (SAD) operating conditions. The results are reported in Table~\ref{table3:cost}. We also report results using SDR, SOCR, and CHR to show the advantages offered by TCR and QCR.

The MATPOWER-solver MIPS~7.0 was used to solve the original ACOPF problem for all test cases. MIPS numerically failed to solve {\tt case588\_sdet\_\_api} (marked with~``*'' in Table~\ref{table3:cost}), we then used the solver FMINCON. Because of invalid limits on $q_{Gg}$ for dispatchable loads modeled as generators with $\pu_{Gg} = 0$, the MATPOWER function {\tt runopf} could not run the {\tt case89\_pegase\_\_api}, {\tt case240\_pserc\_\_api}, and {\tt case1354\_pegase\_\_api} instances (marked with~``**'' in Table~\ref{table3:cost}). We assigned $\pu_{Gg}$ to $-5$~MW for dispatchable loads of these instances.

We solved all the relaxations in MATLAB using CVX~2.2 with the solver MOSEK~9.1.9 and default precision. All the computations were carried out on an Intel(R) Core(TM) i7-3540M CPU @ 3.00 GHz platform. MOSEK ran out of memory when solving SDR for the 1,354-bus instances.

We denote $\lb{\upsilon}$ the best lower bound which is the maximum value among $\hat{\upsilon}_{SOCR}$, $\hat{\upsilon}_{QCR}$, $\hat{\upsilon}_{TCR}$, $\hat{\upsilon}_{CHR}$, and $\hat{\upsilon}_{SDR}$, respective optimal values of SOCR, QCR, TCR, CHR, and SDR. The optimality gap is measured as $100(1-\hat{\upsilon}_R/\ub{\upsilon})$, where $\ub{\upsilon}$ is the upper bound provided by MIPS. The results in Table~\ref{table3:cost} support the following key points:
\begin{enumerate}
	\item TCR is stronger than QCR under TYP conditions. For example, TCR reduces the optimality gap of QCR from 18.80\% to 0.00\% for \texttt{case30\_ieee}.
	
	In each TYP instance, $\ub{\delta}_{\ell} = \pi/6$ for all $\ell \in \Branch$ in~\eqref{eq3:phaselim}. As suggested in~\cite{coffrin2016a}, when phase angle difference bounds are large, QCR is similar to SOCR, thus weaker than TCR.
	
	\item TCR dominates QCR on all but two API instances: \texttt{case3\_lmbd\_\_api} and \texttt{case179\_goc\_\_api}. We note that QCR dominates SDR on the \texttt{case3\_lmbd\_\_api} instance.
	
	From~\eqref{eq3:sf}--\eqref{eq3:st}, $p_{f\ell} \approx \frac{\angle \cplx{v}_k - \angle \cplx{v}_m}{x_\ell} \approx -p_{t\ell}$ for all $\ell = (k,m) \in \Branch$ in a transmission network. By increasing the active power demand, the bound $\ub{\delta}_{\ell} = \pi/6$, $\ell \in \Branch$, might not be as large in API instances as in TYP instances. This could explain why QCR dominates TCR or SDR on some API instances.
	
	\item TCR dominates QCR on 12 out of 20 SAD instances. We note that TCR reduces the optimality gap of QCR from 21.50\% to 0.12\% for \texttt{case14\_ieee\_\_sad}. Although QCR dominates SDR on 3 instances, SDR on average outperforms QCR.
	
	A SAD instance is a TYP instance where $\ub{\delta}_{\ell} = \hat{\delta} \ll \pi/6$ for all $\ell \in \Branch$. The increase in the number of instances where QCR dominates TCR or SDR is explained by the fact that the strength of QCR is sensitive to $\hat{\delta}$.
\end{enumerate}

The computation times reported by MOSEK are also shown in Table~\ref{table3:cost}. The time CVX took to pre-compile a model is not included. We note that TCR is slightly slower than QCR.

\begin{table*}
	\footnotesize
	\centering
	\caption{TCR vs QCR.}\label{table3:cost}
	\resizebox{\linewidth}{!}{
		\begin{tabular*}{1.03\textwidth}{@{}lrr|rrrrr|rrrrr@{}}
			\toprule
			Test case & $\ub{\upsilon}$ [\$/h] & $\lb{\upsilon}$ [\$/h] & \multicolumn{5}{c|}{Optimality gap [\%]} & \multicolumn{5}{c}{Computation time [s]} \\
			\cmidrule(r){4-8} \cmidrule(l){9-13} & 	&  & SOCR &  QCR & TCR & CHR & SDR & SOCR &  QCR & TCR & CHR & SDR \\
			\midrule
			\multicolumn{13}{c}{\it Typical (TYP) instances}\\
			\midrule
			{\tt case3\_lmbd}	&	5 812.64	&	5 789.91	&	1.32	&	1.24	&	0.74	&	0.39	&	0.39	&	0.23	&	0.23	&	0.25	&	0.20	&	0.20	\\
			{\tt case5\_pjm}	&	17 551.89	&	16 635.78	&	14.54	&	14.54	&	12.75	&	5.22	&	5.22	&	0.27	&	0.27	&	0.20	&	0.20	&	0.22	\\
			{\tt case14\_ieee}	&	2 178.08	&	2178.08	&	0.11	&	0.11	&	0.00	&	0.00	&	0.00	&	0.27	&	0.25	&	0.22	&	0.22	&	0.25	\\
			{\tt case24\_ieee\_rts}	&	63 352.20	&	63 352.20	&	0.01	&	0.01	&	0.00	&	0.00	&	0.00	&	0.27	&	0.19	&	0.25	&	0.23	&	0.33	\\
			{\tt case30\_as}	&	803.13	&	803.13	&	0.06	&	0.06	&	0.00	&	0.00	&	0.00	&	0.20	&	0.27	&	0.28	&	0.27	&	0.31	\\
			{\tt case30\_fsr}	&	575.77	&	575.77	&	0.39	&	0.39	&	0.04	&	0.00	&	0.00	&	0.23	&	0.25	&	0.28	&	0.27	&	0.33	\\
			{\tt case30\_ieee}	&	8 208.52	&	8 208.51	&	18.84	&	18.80	&	0.00	&	0.00	&	0.00	&	0.25	&	0.25	&	0.30	&	0.27	&	0.34	\\
			{\tt case39\_epri}	&	138 415.56	&	138 407.21	&	0.55	&	0.54	&	0.20	&	0.01	&	0.01	&	0.34	&	0.36	&	0.41	&	0.38	&	0.52	\\
			{\tt case57\_ieee}	&	37 589.34	&	37 588.31	&	0.16	&	0.16	&	0.01	&	0.00	&	0.00	&	0.25	&	0.30	&	0.45	&	0.41	&	1.08	\\
			{\tt case73\_ieee\_rts}	&	189 764.09	&	189 764.08	&	0.03	&	0.03	&	0.00	&	0.00	&	0.00	&	0.25	&	0.34	&	0.58	&	0.45	&	1.64	\\
			{\tt case89\_pegase}	&	107 285.67	&	106 968.44	&	0.75	&	0.74	&	0.55	&	0.30	&	0.30	&	0.52	&	1.11	&	1.42	&	2.72	&	6.59	\\
			{\tt case118\_ieee}	&	97 213.61	&	97 143.74	&	0.90	&	0.79	&	0.22	&	0.07	&	0.07	&	0.41	&	0.72	&	0.91	&	0.81	&	6.33	\\
			{\tt case162\_ieee\_dtc}	&	108 075.65	&	106 156.53	&	5.94	&	5.91	&	4.98	&	1.78	&	1.78	&	0.56	&	1.14	&	1.95	&	12.25	&	20.14	\\
			{\tt case179\_goc}	&	754 266.42	&	753 724.66	&	0.16	&	0.15	&	0.15	&	0.07	&	0.07	&	0.41	&	0.63	&	1.28	&	1.08	&	21.02	\\
			{\tt case200\_tamu}	&	27 557.57	&	27 557.56	&	0.00	&	0.00	&	0.00	&	0.00	&	0.00	&	0.23	&	0.47	&	0.78	&	1.08	&	14.30	\\
			{\tt case240\_pserc}	&	3 329 670.11	&	3 282 077.65	&	2.77	&	2.73	&	2.60	&	1.43	&	1.43	&	0.86	&	1.72	&	2.81	&	1.66	&	72.06	\\
			{\tt case300\_ieee}	&	565 219.99	&	564 543.99	&	2.62	&	2.59	&	1.17	&	0.12	&	0.12	&	0.94	&	2.00	&	3.19	&	3.39	&	109.63	\\
			{\tt case500\_tamu}	&	72 578.30	&	71 048.42	&	5.38	&	5.38	&	4.39	&	2.11	&	2.11	&	1.28	&	3.58	&	4.94	&	5.19	&	445.89	\\
			{\tt case588\_sdet}	&	313 139.78	&	310 507.73	&	2.18	&	2.14	&	1.61	&	0.84	&	0.85	&	1.55	&	3.34	&	5.03	&	7.01	&	859.56	\\
			{\tt case1354\_pegase}	&	1 258 844.00	&	1 251 841.24	&	1.57	&	1.57	&	1.23	&	0.56	&	--	&	9.70	&	23.48	&	25.22	&	34.06	&	--	\\
			{\bf Average}	&	&	&	\bf 2.33	&	\bf 2.31	&	\bf 1.68	&	\bf 0.75	&	\bf --	&	\bf 3.86	&	\bf 9.23	&	\bf 10.47	&	\bf 14.15	&	\bf --	\\
			\midrule
			\multicolumn{13}{c}{\it Active power increase (API) instances}\\
			\midrule
			{\tt case3\_lmbd\_\_api}	&	11 242.13	&	10 450.57	&	9.32	&	7.04	&	7.90	&	7.34	&	7.34	&	0.25	&	0.25	&	0.19	&	0.17	&	0.19	\\
			{\tt case5\_pjm\_\_api}	&	76 377.42	&	76 182.35	&	4.09	&	4.09	&	3.22	&	0.26	&	0.26	&	0.23	&	0.42	&	0.25	&	0.23	&	0.30	\\
			{\tt case14\_ieee\_\_api}	&	5 999.36	&	5 999.36	&	5.13	&	5.13	&	0.57	&	0.00	&	0.00	&	0.30	&	0.25	&	0.25	&	0.28	&	0.23	\\
			{\tt case24\_ieee\_rts\_\_api}	&	134 948.17	&	132 161.97	&	17.87	&	13.04	&	6.01	&	2.06	&	2.06	&	0.28	&	0.34	&	0.39	&	0.27	&	0.36	\\
			{\tt case30\_as\_\_api}	&	4 996.21	&	4 925.85	&	44.60	&	44.60	&	42.34	&	1.41	&	1.41	&	0.22	&	0.31	&	0.42	&	0.38	&	0.59	\\
			{\tt case30\_fsr\_\_api}	&	701.15	&	699.16	&	2.76	&	2.75	&	2.01	&	0.28	&	0.28	&	0.27	&	0.38	&	0.38	&	0.31	&	0.50	\\
			{\tt case30\_ieee\_\_api}	&	18 043.92	&	18 043.87	&	5.45	&	5.45	&	0.36	&	0.00	&	0.00	&	0.48	&	0.31	&	0.33	&	0.31	&	0.44	\\
			{\tt case39\_epri\_\_api}	&	249 747.58	&	249 281.70	&	1.73	&	1.71	&	1.05	&	0.19	&	0.19	&	0.33	&	0.34	&	0.39	&	0.39	&	0.69	\\
			{\tt case57\_ieee\_\_api}	&	49 296.69	&	49 294.98	&	0.08	&	0.08	&	0.03	&	0.00	&	0.00	&	0.23	&	0.31	&	0.44	&	0.45	&	1.01	\\
			{\tt case73\_ieee\_rts\_\_api}	&	422 726.27	&	410 415.90	&	12.88	&	11.84	&	9.34	&	2.91	&	2.91	&	0.25	&	0.53	&	0.58	&	0.66	&	2.09	\\
			{\tt case89\_pegase\_\_api}	&	**134 276.98	&	118 382.04	&	13.40	&	13.39	&	13.04	&	11.84	&	11.84	&	0.59	&	1.08	&	1.78	&	3.14	&	8.05	\\
			{\tt case118\_ieee\_\_api}	&	242 054.01	&	215 037.82	&	28.81	&	28.70	&	26.41	&	11.16	&	11.16	&	0.38	&	0.70	&	0.94	&	0.91	&	7.61	\\
			{\tt case162\_ieee\_dtc\_\_api}	&	120 996.09	&	119 276.13	&	4.36	&	4.35	&	3.65	&	1.42	&	1.42	&	0.63	&	1.17	&	1.83	&	12.47	&	22.56	\\
			{\tt case179\_goc\_\_api}	&	1 932 120.36	&	1 921 462.02	&	9.88	&	5.92	&	6.78	&	0.55	&	0.55	&	0.56	&	1.11	&	1.89	&	1.31	&	22.41	\\
			{\tt case200\_tamu\_\_api}	&	36 763.28	&	36 763.26	&	0.02	&	0.02	&	0.00	&	0.00	&	0.00	&	0.50	&	0.89	&	1.53	&	1.80	&	24.86	\\
			{\tt case240\_pserc\_\_api}	&	**4 768 555.29	&	4 754 073.18	&	0.73	&	0.69	&	0.62	&	0.30	&	0.30	&	0.84	&	1.98	&	3.09	&	2.20	&	70.48	\\
			{\tt case300\_ieee\_\_api}	&	650 147.23	&	649 628.70	&	0.89	&	0.83	&	0.36	&	0.08	&	0.08	&	0.95	&	2.28	&	3.30	&	25.92	&	107.69	\\
			{\tt case500\_tamu\_\_api}	&	42 775.62	&	42 288.76	&	2.91	&	2.91	&	2.15	&	1.14	&	1.14	&	1.27	&	3.39	&	4.27	&	5.08	&	424.30	\\
			{\tt case588\_sdet\_\_api}	&	*394 758.02	&	393 096.82	&	1.65	&	1.61	&	0.99	&	0.42	&	0.42	&	1.89	&	4.55	&	6.75	&	8.25	&	876.42	\\
			{\tt case1354\_pegase\_\_api}	&	**1 487 103.11	&	1 482 481.21	&	0.66	&	0.66	&	0.47	&	0.31	&	--	&	11.58	&	26.53	&	25.38	&	33.66	&	--	\\
			{\bf Average}	&	&	&	\bf 3.46	&	\bf 3.22	&	\bf 2.70	&	\bf 1.08	&	\bf --	&	\bf 4.56	&	\bf 10.49	&	\bf 10.79	&	\bf 15.96	&	\bf --	\\
			\midrule
			\multicolumn{13}{c}{\it Small angle difference (SAD) instances}\\
			\midrule
			{\tt case3\_lmbd\_\_sad}	&	5 959.33	&	5 874.07	&	3.74	&	1.43	&	2.42	&	1.86	&	1.86	&	0.56	&	0.38	&	0.69	&	0.41	&	0.39	\\
			{\tt case5\_pjm\_\_sad}	&	26 115.20	&	26 115.20	&	3.62	&	0.99	&	3.28	&	0.00	&	0.00	&	0.80	&	0.55	&	0.52	&	0.56	&	0.52	\\
			{\tt case14\_ieee\_\_sad}	&	2 777.30	&	2 774.79	&	21.54	&	21.49	&	0.12	&	0.09	&	0.09	&	0.34	&	0.33	&	0.38	&	0.44	&	0.30	\\
			{\tt case24\_ieee\_rts\_\_sad}	&	76 943.25	&	74 691.09	&	9.55	&	2.93	&	6.93	&	4.36	&	4.36	&	0.91	&	0.59	&	0.88	&	0.69	&	0.97	\\
			{\tt case30\_as\_\_sad}	&	897.49	&	895.34	&	7.97	&	2.31	&	0.43	&	0.24	&	0.24	&	0.44	&	0.44	&	0.47	&	0.44	&	0.61	\\
			{\tt case30\_fsr\_\_sad}	&	576.79	&	576.68	&	0.47	&	0.41	&	0.11	&	0.02	&	0.02	&	0.38	&	0.55	&	0.48	&	0.45	&	0.61	\\
			{\tt case30\_ieee\_\_sad}	&	8 208.52	&	8 208.51	&	9.69	&	5.93	&	0.00	&	0.00	&	0.00	&	0.41	&	0.45	&	0.55	&	0.52	&	0.61	\\
			{\tt case39\_epri\_\_sad}	&	148 354.42	&	148 324.17	&	0.66	&	0.21	&	0.09	&	0.02	&	0.02	&	0.75	&	0.58	&	0.67	&	0.73	&	1.00	\\
			{\tt case57\_ieee\_\_sad}	&	38 663.88	&	38 646.23	&	0.70	&	0.35	&	0.18	&	0.05	&	0.05	&	0.73	&	0.81	&	1.27	&	1.30	&	3.19	\\
			{\tt case73\_ieee\_rts\_\_sad}	&	227 745.75	&	221 967.33	&	6.75	&	2.54	&	5.54	&	2.75	&	2.75	&	0.78	&	1.09	&	1.55	&	1.49	&	6.13	\\
			{\tt case89\_pegase\_\_sad}	&	107 285.67	&	106 968.65	&	0.74	&	0.70	&	0.55	&	0.30	&	0.30	&	1.06	&	2.30	&	3.23	&	5.59	&	12.33	\\
			{\tt case118\_ieee\_\_sad}	&	105 216.69	&	101 813.57	&	8.25	&	6.83	&	6.88	&	3.23	&	3.23	&	0.80	&	1.05	&	2.25	&	2.05	&	15.03	\\
			{\tt case162\_ieee\_dtc\_\_sad}	&	108 695.95	&	106 437.44	&	6.48	&	6.35	&	5.51	&	2.08	&	2.08	&	0.53	&	1.51	&	2.63	&	17.27	&	31.09	\\
			{\tt case179\_goc\_\_sad}	&	762 541.29	&	755 322.23	&	1.12	&	1.00	&	1.11	&	0.95	&	0.95	&	0.53	&	1.00	&	1.77	&	1.36	&	26.67	\\
			{\tt case200\_tamu\_\_sad}	&	27 557.57	&	27 557.56	&	0.00	&	0.00	&	0.00	&	0.00	&	0.00	&	0.41	&	0.78	&	1.20	&	1.66	&	20.34	\\
			{\tt case240\_pserc\_\_sad}	&	3 407 087.72	&	3 289 232.26	&	4.98	&	4.41	&	4.79	&	3.46	&	3.46	&	0.84	&	1.49	&	2.89	&	1.81	&	72.39	\\
			{\tt case300\_ieee\_\_sad}	&	565 712.85	&	564 913.60	&	2.67	&	2.46	&	1.26	&	0.14	&	0.14	&	0.92	&	1.73	&	3.16	&	3.41	&	108.97	\\
			{\tt case500\_tamu\_\_sad}	&	79 233.96	&	73 222.36	&	7.91	&	7.89	&	7.74	&	7.59	&	7.59	&	1.55	&	3.72	&	4.91	&	5.72	&	449.78	\\
			{\tt case588\_sdet\_\_sad}	&	329 860.72	&	311 564.46	&	6.94	&	6.24	&	6.34	&	5.55	&	5.55	&	1.47	&	3.50	&	5.31	&	7.39	&	683.01	\\
			{\tt case1354\_pegase\_\_sad}	&	1 258 848.13	&	1 251 837.86	&	1.57	&	1.55	&	1.23	&	0.56	&	--	&	16.03	&	27.17	&	24.14	&	34.61	&	--	\\
			{\bf Average}	&	&	&	\bf 3.99	&	\bf 3.58	&	\bf 3.31	&	\bf 2.46	&	\bf --	&	\bf 6.06	&	\bf 10.58	&	\bf 10.34	&	\bf 14.85	&	\bf --	\\
			\bottomrule
		\end{tabular*}
	}
\end{table*}

\section{Conclusion}\label{sec3:conclusion}
Using PGLib-OPF test cases with up to 1,354 buses, we show that on average, TCR dominates QCR. Furthermore, for test cases under typical operating conditions, TCR is stronger than QCR. Overall, both in optimality gap and computation time, TCR is a better trade-off than QCR.


\end{document}